\documentclass[11pt]{article}

\usepackage{amsmath,amsthm,amscd,amsfonts,amssymb}
\usepackage{graphicx,slashbox}
\usepackage{subfigure}

\usepackage{a4wide}

\newtheorem{thm}{Theorem}[section]
\newtheorem{definition}[thm]{Definition}
\newtheorem{remark}[thm]{Remark}

\theoremstyle{definition}
\newtheorem{example}[thm]{Example}

\def\a{\alpha}

\def\LI{{_aI_t^\a}}
\def\RI{{_tI_b^\a}}

\def\LIz{{_0I_t^{0.5}}}
\def\LIza{{_0I_t^{\a}}}

\newenvironment{keywords}{\begin{center}
\begin{minipage}[c]{13.4cm} {\bf Keywords:}} {\end{minipage}
\end{center}}

\newenvironment{msc}{\begin{center}
\begin{minipage}[c]{13.4cm} {\bf MSC 2010:}} {\end{minipage}
\end{center}}

\begin{document}

\title{Approximation of fractional integrals
by means of derivatives\thanks{Submitted 15-10-2011; revised 24-01-2012;
accepted 25-01-2012; for publication in 
\emph{Computers and Mathematics with Applications}.
Part of the first author's Ph.D.,
which is carried out at the University of Aveiro under the
\emph{Doctoral Program in Mathematics and Applications} (PDMA)
of Universities of Aveiro and Minho.}}

\author{Shakoor Pooseh\\
\texttt{spooseh@ua.pt}
\and Ricardo Almeida\\
\texttt{ricardo.almeida@ua.pt}
\and Delfim F. M. Torres\\
\texttt{delfim@ua.pt}}

\date{Center for Research and Development in Mathematics and Applications\\
Department of Mathematics, University of Aveiro, 3810-193 Aveiro, Portugal}

\maketitle


\begin{abstract}
We obtain a new decomposition of the Riemann--Liouville operators
of fractional integration as a series involving derivatives (of integer order).
The new formulas are valid for functions of class $C^n$,
$n \in \mathbb{N}$, and allow us to develop suitable numerical
approximations with known estimations for the error.
The usefulness of the obtained results,
in solving fractional integral equations
and fractional problems of the calculus of variations,
is illustrated.
\end{abstract}

\begin{msc}
26A33, 33F05.
\end{msc}



\section{Introduction}

Let ${_aI_t^1}x(t) := \int_a^t x(\tau) d\tau$. If the equation
\begin{equation}
\label{eq:nfoldI}
{_aI_t^n}x(t) = \frac{1}{(n-1)!}\int_a^t (t-\tau)^{n-1}x(\tau) d\tau
\end{equation}
is true for the $n$-fold integral, $n \in \mathbb{N}$, then
\begin{equation*}
\begin{split}
{_aI_t^{n+1}}x(t)
&= {_aI_t^1} \left(\frac{1}{(n-1)!}\int_a^t (t-\tau)^{n-1}x(\tau) d\tau\right)\\
&= \int_a^t \left(\frac{1}{(n-1)!}\int_a^\xi (\xi-\tau)^{n-1}x(\tau) d\tau\right)d\xi.
\end{split}
\end{equation*}
Interchanging the order of integration gives
\begin{equation*}
{_aI_t^{n+1}}x(t) = \frac{1}{n!}\int_a^t (t-\tau)^{n}x(\tau) d\tau .
\end{equation*}
Since, by definition, \eqref{eq:nfoldI} is true for $n=1$,
so it is also true for all $n \in \mathbb{N}$ by induction.
The (left Riemann--Liouville) fractional integral of $x(t)$ of order
$\alpha > 0$ is then naturally defined, as an extension of \eqref{eq:nfoldI},
with the help of Euler's Gamma function $\Gamma$:
\begin{equation}
\label{eq:nfoldI:alpha}
{_aI_t^{\alpha}}x(t) := \frac{1}{\Gamma(\alpha)}\int_a^t (t-\tau)^{\alpha-1}x(\tau) d\tau .
\end{equation}
The study of fractional integrals \eqref{eq:nfoldI:alpha}
is a two hundred years old subject that is part of
a branch of mathematical analysis called \emph{Fractional Calculus}
\cite{Kilbas,Miller,book:Samko}. Recently, due to its many applications
in science and engineering, there has been an increase of interest
in the study of fractional calculus \cite{TM}.
Fractional integrals appear naturally in many different contexts,
\textrm{e.g.}, when dealing with fractional variational problems
or fractional optimal control \cite{VibString,APT,FT,MT,Dorota}.
As is frequently observed, solving such equations
analytically can be a difficult task,
even impossible in some cases. One way to overcome
the problem consists to apply numerical methods, \textrm{e.g.},
using Riemann sums to approximate the fractional operators.
We refer the reader to \cite{AlD1,Atan2,Li,Yang} and references therein.

Here we obtain a simple and effective approximation for
fractional integrals. The paper is organized as follows.
First, in Section~\ref{secpre}, we fix some notation by recalling
the basic definitions of fractional calculus.
In Section~\ref{secDecomp} we obtain a decomposition formula for the left
and right fractional integrals of functions of class $C^n$ (Theorems~\ref{TheoDecomp}
and \ref{TheoDecomp2}). The error derived by these approximations is studied
in Section~\ref{secError}. In Section~\ref{secExp} we consider several examples,
where we determine the exact expression of the fractional integrals for some functions,
and compare them with numerical approximations of different types.
We end with Section~\ref{secExp2} of applications, where we solve numerically,
by means of the obtained approximations, an equation depending on a fractional integral;
and a fractional problem of the calculus of variations.


\section{Preliminaries}
\label{secpre}

We fix notations by recalling the  basic concepts
(see, \textrm{e.g.}, \cite{Kilbas}).

\begin{definition}
Let $x(\cdot)$ be an integrable function in $[a,b]$ and $\a>0$.
The left Riemann--Liouville fractional integral of order $\alpha$ is given by
$$
\LI x(t)=\frac{1}{\Gamma(\alpha)}
\int_a^t (t-\tau)^{\alpha-1}x(\tau)d\tau,
\quad t\in [a,b],
$$
while the right Riemann--Liouville fractional integral
of order $\alpha$ is given by
$$
\RI x(t)=\frac{1}{\Gamma(\alpha)}\int_t^b (\tau-t)^{\a-1}x(\tau)d\tau,
\quad t\in [a,b].
$$
\end{definition}

If $\a>0$, $\beta>0$, and $x\in L_p(a,b)$, $1\leq p \leq \infty$, then
$$
\LI {_aI_t^\beta}x(t)={_aI_t^{\a+\beta}}x(t)\ \mbox{ and }
\ \RI {_tI_b^\beta}x(t)={_tI_b^{\a+\beta}}x(t)
$$
almost everywhere. The equalities
hold for all $t\in[a,b]$ if in addition $\a+\beta>1$.


\section{A decomposition for the fractional integral}
\label{secDecomp}

For analytical functions, we can rewrite a fractional integral
as a series involving integer derivatives only.
If $x$ is analytic in $[a,b]$, then
\begin{equation}
\label{analytical}
\LI x(t)=\frac{1}{\Gamma(\a)}
\sum_{k=0}^\infty\frac{(-1)^k(t-a)^{k+\a}}{(k+\a)k!}x^{(k)}(t)
\end{equation}
for all $t\in[a,b]$ (\textrm{cf.} Eq. (3.44) in \cite{Miller}).
From the numerical point of view, one considers finite sums
and the following approximation:
\begin{equation}
\label{analytical2}
\LI x(t)\approx\frac{1}{\Gamma(\a)}\sum_{k=0}^N
\frac{(-1)^k(t-a)^{k+\a}}{(k+\a)k!}x^{(k)}(t).
\end{equation}

One problem with formula \eqref{analytical}
is the restricted class of functions
where it is valid. In applications,
this approach may not be suitable. The main aim of this paper
is to present a new decomposition formula for functions of class $C^n$.
Before we give the result in its full extension, we explain the method for $n=3$.
To that purpose, let $x\in C^3[a,b]$. Using integration by parts three times, we deduce that
$$
\LI x(t)=\frac{(t-a)^\a}{\Gamma(\a+1)}x(a)
+\frac{(t-a)^{\a+1}}{\Gamma(\a+2)}x'(a)
+\frac{(t-a)^{\a+2}}{\Gamma(\a+3)}x''(a)
+\frac{1}{\Gamma(\a+3)}\int_a^t (t-\tau)^{\a+2}x^{(3)}(\tau)d\tau.
$$
By the binomial formula, we can rewrite the fractional integral as
\begin{multline*}
\LI x(t)=\frac{(t-a)^\a}{\Gamma(\a+1)}x(a)
+\frac{(t-a)^{\a+1}}{\Gamma(\a+2)}x'(a)
+\frac{(t-a)^{\a+2}}{\Gamma(\a+3)}x''(a)\\
+\frac{(t-a)^{\a+2}}{\Gamma(\a+3)}\sum_{p=0}^\infty
\frac{\Gamma(p-\a-2)}{\Gamma(-\a-2)p!(t-a)^p}
\int_a^t (\tau-a)^px^{(3)}(\tau)d\tau.
\end{multline*}
The rest of the procedure follows the same pattern: decompose the sum into
a first term plus the others, and integrate by parts. Then we obtain
\begin{equation*}
\begin{split}
\LI x(t) &= \frac{(t-a)^\a}{\Gamma(\a+1)}x(a)
+\frac{(t-a)^{\a+1}}{\Gamma(\a+2)}x'(a)+\frac{(t-a)^{\a+2}}{\Gamma(\a+3)}x''(t)
\left[1+\sum_{p=1}^\infty\frac{\Gamma(p-\a-2)}{\Gamma(-\a-2)p!}\right]\\
&\quad +\frac{(t-a)^{\a+2}}{\Gamma(\a+2)}\sum_{p=1}^\infty
\frac{\Gamma(p-\a-2)}{\Gamma(-\a-1)(p-1)!(t-a)^p}
\int_a^t (\tau-a)^{p-1}x''(\tau)d\tau\\
&= \frac{(t-a)^\a}{\Gamma(\a+1)}x(a)
+\frac{(t-a)^{\a+1}}{\Gamma(\a+2)}x'(t)\left[1+\sum_{p=2}^\infty
\frac{\Gamma(p-\a-2)}{\Gamma(-\a-1)(p-1)!}\right]\\
&\quad +\frac{(t-a)^{\a+2}}{\Gamma(\a+3)}x''(t)\left[1+\sum_{p=1}^\infty
\frac{\Gamma(p-\a-2)}{\Gamma(-\a-2)p!}\right]\\
&\quad +\frac{(t-a)^{\a+2}}{\Gamma(\a+1)}\sum_{p=2}^\infty
\frac{\Gamma(p-\a-2)}{\Gamma(-\a)(p-2)!(t-a)^p}
\int_a^t (\tau-a)^{p-2}x'(\tau)d\tau
\end{split}
\end{equation*}
\begin{equation*}
\begin{split}
&= \frac{(t-a)^\a}{\Gamma(\a+1)}x(t)\left[1+\sum_{p=3}^\infty
\frac{\Gamma(p-\a-2)}{\Gamma(-\a)(p-2)!}\right]
+\frac{(t-a)^{\a+1}}{\Gamma(\a+2)}x'(t)\left[1+\sum_{p=2}^\infty
\frac{\Gamma(p-\a-2)}{\Gamma(-\a-1)(p-1)!}\right]\\
&\quad +\frac{(t-a)^{\a+2}}{\Gamma(\a+3)}x''(t)\left[1+\sum_{p=1}^\infty
\frac{\Gamma(p-\a-2)}{\Gamma(-\a-2)p!}\right]\\
&\quad +\frac{(t-a)^{\a+2}}{\Gamma(\a)}\sum_{p=3}^\infty
\frac{\Gamma(p-\a-2)}{\Gamma(-\a+1)(p-3)!(t-a)^p}\int_a^t (\tau-a)^{p-3}x(\tau)d\tau.
\end{split}
\end{equation*}
Therefore, we can expand $\LI x(t)$ as
\begin{multline}
\label{def:n=3}
\LI x(t)=A_0(\alpha)(t-a)^\a x(t)+ A_1(\alpha)(t-a)^{\a+1} x'(t)
+ A_2(\alpha)(t-a)^{\a+2} x''(t)\\
+ \sum_{p=3}^\infty B(\a,p)(t-a)^{\a+2-p}V_p(t),
\end{multline}
where
$$
\begin{array}{ll}
A_0(\a)&=\displaystyle\frac{1}{\Gamma(\a+1)}\left[1+\sum_{p=3}^\infty
\frac{\Gamma(p-\a-2)}{\Gamma(-\a)(p-2)!}\right],\\
A_1(\a)&=\displaystyle\frac{1}{\Gamma(\a+2)}\left[1+\sum_{p=2}^\infty
\frac{\Gamma(p-\a-2)}{\Gamma(-\a-1)(p-1)!}\right],\\
A_2(\a)&=\displaystyle\frac{1}{\Gamma(\a+3)}\left[1+\sum_{p=1}^\infty
\frac{\Gamma(p-\a-2)}{\Gamma(-\a-2)p!}\right],\\
\end{array}
$$
\begin{equation}
\label{def:B3}
B(\a,p)=\frac{\Gamma(p-\a-2)}{\Gamma(\a)\Gamma(1-\a)(p-2)!},
\end{equation}
and
\begin{equation}
\label{def:Vp3}
V_p(t)=\int_a^t (p-2)(\tau-a)^{p-3}x(\tau)d\tau.
\end{equation}

\begin{remark}
Function $V_p$ given by \eqref{def:Vp3} may be defined
as the solution of the differential equation
$$
\left\{
\begin{array}{l}
V_p'(t)=(p-2)(t-a)^{p-3}x(t)\\
V_p(a)=0\\
\end{array}\right.
$$
for $p=3,4,\ldots$
\end{remark}

\begin{remark}
When $\a$ is not an integer, we may use
Euler's reflection formula (\textrm{cf.} \cite{Beals})
$$
\Gamma(\a)\Gamma(1-\a)=\frac{\pi}{\sin(\pi\a)}
$$
to simplify expression $B(\a,p)$ in \eqref{def:B3}.
\end{remark}

Following the same reasoning, we are able to deduce a general formula
of decomposition for fractional integrals, depending on the order
of smoothness of the test function.

\begin{thm}
\label{TheoDecomp}
Let $n\in\mathbb N$ and $x\in C^n[a,b]$. Then
\begin{equation}
\label{ExpanDecomp}
\LI x(t)=\sum_{i=0}^{n-1}A_i(\alpha)(t-a)^{\a+i} x^{(i)}(t)
+\sum_{p=n}^\infty B(\a,p)(t-a)^{\a+n-1-p}V_p(t),
\end{equation}
where
\begin{equation}
\label{def:B}
\begin{array}{ll}
A_i(\a)&=\displaystyle\frac{1}{\Gamma(\a+i+1)}\left[1+\sum_{p=n-i}^\infty
\frac{\Gamma(p-\a-n+1)}{\Gamma(-\a-i)(p-n+1+i)!}\right],
\quad i = 0, \ldots, n-1,\\
B(\a,p)&=\displaystyle\frac{\Gamma(p-\a-n+1)}{\Gamma(\a)\Gamma(1-\a)(p-n+1)!},
\end{array}
\end{equation}
and
\begin{equation}
\label{def:Vp}
V_p(t)=\int_a^t (p-n+1)(\tau-a)^{p-n}x(\tau)d\tau,
\end{equation}
$p = n, n+1, \ldots$
\end{thm}

A remark about the convergence of the series in $A_i(\a)$,
for $i\in\{0,\ldots,n-1\}$, is in order. Since
\begin{equation}
\label{eq:conv:ser}
\begin{array}{ll}
\displaystyle \sum_{p=n-i}^\infty\frac{\Gamma(p-\a-n+1)}{\Gamma(-\a-i)(p-n+1+i)!}
&= \displaystyle\sum_{p=0}^\infty\frac{\Gamma(p-\a-i)}{\Gamma(-\a-i) p!}-1\\
&={_2F_1} (-\a-i,-,-,1)-1,\\
\end{array}
\end{equation}
where ${_2F_1}$ denotes the hypergeometric function,
and because $\a+i>0$, we conclude that \eqref{eq:conv:ser}
converges absolutely (\textrm{cf.} Theorem 2.1.2 in \cite{Andrews}).
In fact, we may use Eq. (2.1.6) in \cite{Andrews} to conclude that
$$
\sum_{p=n-i}^\infty\frac{\Gamma(p-\a-n+1)}{\Gamma(-\a-i)(p-n+1+i)!}=-1.
$$
Therefore, the first $n$ terms of our decomposition \eqref{ExpanDecomp} vanish.
However, because of numerical reasons, we do not follow this procedure here.
Indeed, only finite sums of these coefficients are to be taken,
and we obtain a better accuracy for the approximation
taking them into account (see Figures~\ref{ExpTkA} and \ref{ExpTk2A}).
More precisely, we consider finite sums up to order $N$,
with $N\geq n$. Thus, our approximation will depend on two parameters:
the order of the derivative $n\in\mathbb N$, and the number
of terms taken in the sum, which is given by $N$.
The left fractional integral is then approximated by
\begin{equation}
\label{Approx:LI}
\LI x(t)\approx \sum_{i=0}^{n-1}A_i(\alpha,N)(t-a)^{\a+i} x^{(i)}(t)
+\sum_{p=n}^N B(\a,p)(t-a)^{\a+n-1-p}V_p(t),
\end{equation}
where
\begin{equation}
\label{Def:A}
A_i(\a,N)=\frac{1}{\Gamma(\a+i+1)}\left[1+\sum_{p=n-i}^N
\frac{\Gamma(p-\a-n+1)}{\Gamma(-\a-i)(p-n+1+i)!}\right],
\end{equation}
and $B(\a,p)$ and $V_p(t)$ are given by \eqref{def:B}
and \eqref{def:Vp}, respectively.

To measure the errors made by neglecting
the remaining terms, observe that
\begin{equation}
\label{error:A}
\begin{split}
\frac{1}{\Gamma(\a+i+1)}&\sum_{p=N+1}^\infty\frac{\Gamma(p-\a-n+1)}{\Gamma(-\a-i)(p-n+1+i)!}
=\frac{1}{\Gamma(\a+i+1)}\sum_{p=N-n+2+i}^\infty\frac{\Gamma(p-\a-i)}{\Gamma(-\a-i)p!}\\
&= \frac{1}{\Gamma(\a+i+1)}\left[  {_2F_1} (-\a-i,-,-,1)
- \sum_{p=0}^{N-n+1+i}\frac{\Gamma(p-\a-i)}{\Gamma(-\a-i)p!} \right]\\
&=\frac{-1}{\Gamma(\a+i+1)}\sum_{p=0}^{N-n+i+1}\frac{\Gamma(p-\a-i)}{\Gamma(-\a-i)p!}.
\end{split}
\end{equation}
Similarly,
\begin{equation}
\label{error:B}
\frac{1}{\Gamma(\a)\Gamma(1-\a)}\sum_{p=N+1}^\infty\frac{\Gamma(p-\a-n+1)}{(p-n+1)!}
=\frac{-1}{\Gamma(\a)\Gamma(1-\a)}\sum_{p=0}^{N-n+1}\frac{\Gamma(p-\a)}{p!}.
\end{equation}
In Tables~\ref{tab1} and \ref{tab2} we exemplify some values
for \eqref{error:A} and \eqref{error:B}, respectively, with $\a=0.5$
and for different values of $N$, $n$ and $i$. Observe that the errors
only depend on the values of $N-n$ and $i$ for \eqref{error:A},
and on the value of  $N-n$ for \eqref{error:B}.


\begin{table}[!ht]
\center
\begin{tabular}{|c|c|c|c|c|c|}
\hline
\backslashbox{$i$}{$N-n$} & 0 & 1 & 2 & 3 & 4  \\
 \hline
0 & -0.5642 & -0.4231 & -0.3526 & -0.3085 & -0.2777  \\
 \hline
1 &0.09403 & 0.04702 & 0.02938 & 0.02057 & 0.01543 \\
 \hline
2 & -0.01881 & -0.007052& -0.003526& -0.002057& -0.001322\\
 \hline
3 & 0.003358& 0.001007 & 0.0004198 & 0.0002099 & 0.0001181 \\
 \hline
4 & -0.0005224 & -0.0001306 & -0.00004664& -0.00002041 & -0.00001020\\
 \hline
5 & $7.124 \times10^{-5}$& $1.526\times10^{-5}$& $4.770\times10^{-6}$
  &$ 1.855\times10^{-6}$ & $8.347\times10^{-7}$ \\
 \hline
\end{tabular}
\caption{Values of error \eqref{error:A} for $\a=0.5$.}
\label{tab1}
\end{table}


\begin{table}[!ht]
\center
\begin{tabular}{|c|c|c|c|c|c|}
\hline
 $N-n$ & 0 & 1 & 2 & 3 & 4  \\
 \hline
 & 0.5642 & 0.4231 & 0.3526 & 0.3085 & 0.2777  \\
 \hline
\end{tabular}
\caption{Values of error \eqref{error:B} for $\a=0.5$.}
\label{tab2}
\end{table}


Everything done so far is easily adapted
to the right fractional integral.
In particular, one has:

\begin{thm}
\label{TheoDecomp2}
Let $n\in\mathbb N$ and $x\in C^n[a,b]$. Then
$$
\RI x(t)=\sum_{i=0}^{n-1}A_i(\alpha)(b-t)^{\a+i} x^{(i)}(t)
+\sum_{p=n}^\infty B(\a,p)(b-t)^{\a+n-1-p}W_p(t),
$$
where
\begin{equation*}
\begin{split}
A_i(\a)&=\frac{(-1)^i}{\Gamma(\a+i+1)}\left[1+\sum_{p=n-i}^\infty
\frac{\Gamma(p-\a-n+1)}{\Gamma(-\a-i)(p-n+1+i)!}\right],\\
B(\a,p)&=\frac{(-1)^n\Gamma(p-\a-n+1)}{\Gamma(\a)\Gamma(1-\a)(p-n+1)!},
\qquad W_p(t) = \int_t^b (p-n+1)(b-\tau)^{p-n}x(\tau)d\tau.
\end{split}
\end{equation*}
\end{thm}


\section{Error analysis}
\label{secError}

In the previous section we deduced an approximation formula for
the left fractional integral (Eq.~\eqref{Approx:LI}).
The order of magnitude of the coefficients that we ignore during
this procedure is small for the examples that we have chosen
(Tables~\ref{tab1} and \ref{tab2}). The aim of this section
is to obtain an estimation for the error, when considering sums
up to order $N$. We proved that
\begin{multline*}
\LI x(t)=\frac{(t-a)^\a}{\Gamma(\a+1)}x(a)
+\cdots+\frac{(t-a)^{\a+n-1}}{\Gamma(\a+n)}x^{(n-1)}(a)\\
+\frac{(t-a)^{\a+n-1}}{\Gamma(\a+n)}
\int_a^t \left(1-\frac{\tau-a}{t-a}\right)^{\a+n-1}x^{(n)}(\tau)d\tau.
\end{multline*}
Expanding up to order $N$ the binomial, we get
$$
\left(1-\frac{\tau-a}{t-a}\right)^{\a+n-1}
=\sum_{p=0}^N\frac{\Gamma(p-\a-n+1)}{\Gamma(1-\a-n)\,p!}\left(\frac{\tau
-a}{t-a}\right)^p+R_N(\tau),
$$
where
$$
R_N(\tau)=\sum_{p=N+1}^\infty\frac{\Gamma(p-\a-n+1)}{\Gamma(1
-\a-n)\,p!}\left(\frac{\tau-a}{t-a}\right)^p.
$$
Since $\tau\in[a,t]$, we easily deduce an upper bound for $R_N(\tau)$:
\begin{equation*}
\begin{split}
\left|R_N(\tau)\right| &\leq  \sum_{p=N+1}^\infty\left|
\frac{\Gamma(p-\a-n+1)}{\Gamma(1-\a-n)\,p!} \right|
= \sum_{p=N+1}^\infty\left| \binom{\a+n-1}{p} \right|
\leq  \sum_{p=N+1}^\infty \frac{e^{(\a+n-1)^2+\a+n-1}}{p^{\a+n}}\\
&\leq \int_{N}^\infty \frac{e^{(\a+n-1)^2+\a+n-1}}{p^{\a+n}}dp
=\frac{e^{(\a+n-1)^2+\a+n-1}}{(\a+n-1)N^{\a+n-1}}.
\end{split}
\end{equation*}
Thus, we obtain an estimation for the error $E_{tr}(\cdot)$:
$$
\left| E_{tr}(t)\right|\leq L_n \frac{(t-a)^{\a+n}
e^{(\a+n-1)^2+\a+n-1}}{\Gamma(\a+n)(\a+n-1)N^{\a+n-1}},
$$
where $\displaystyle L_n=\max_{\tau\in[a,t]}|x^{(n)}(\tau)|$.


\section{Numerical examples}
\label{secExp}

\begin{figure}[!ht]
  \begin{center}
    \subfigure[$\LIz(t^3)$]{\label{ExpTk}\includegraphics[scale=0.5]{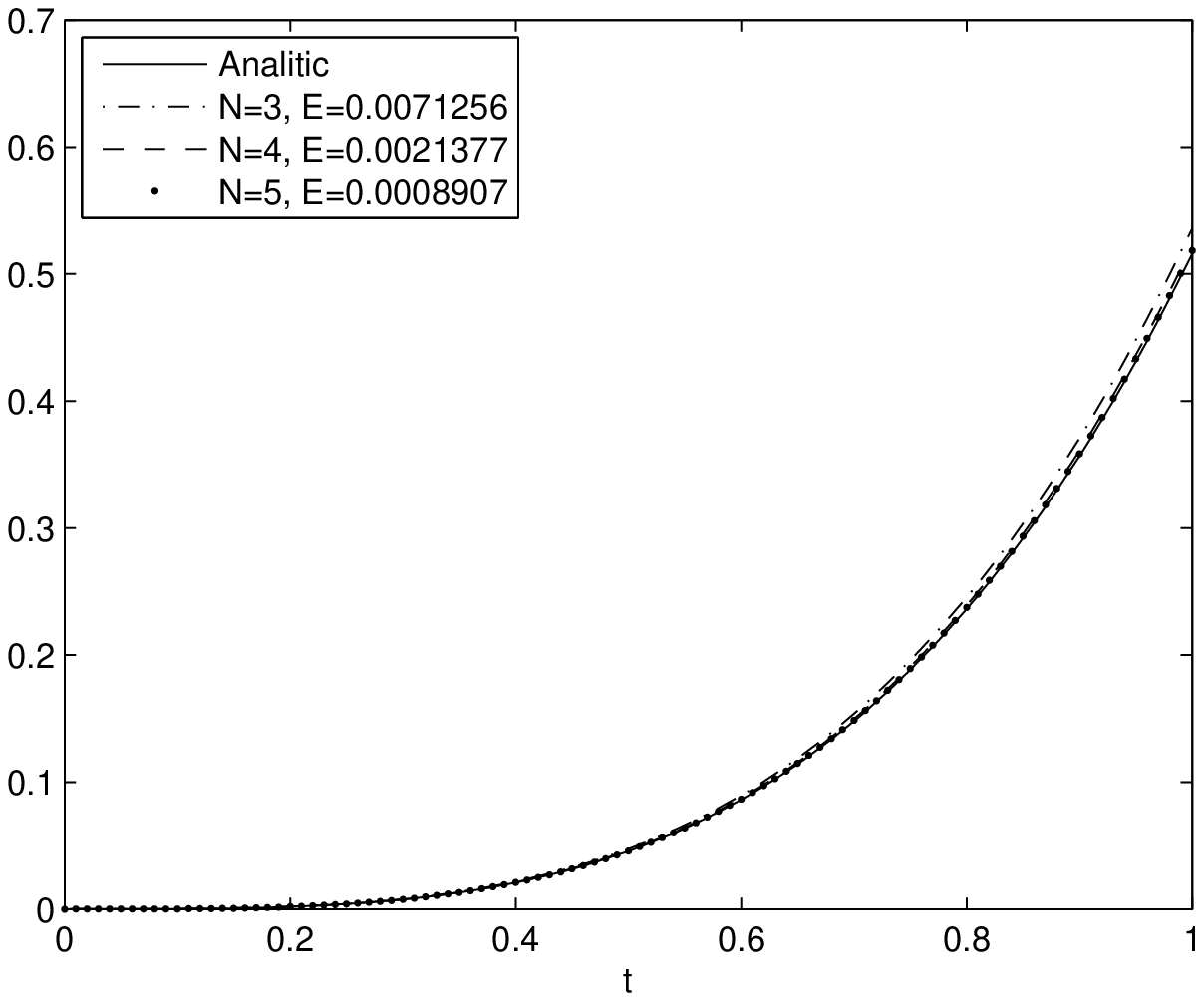}}
    \subfigure[$\LIz(t^{10})$]{\label{ExpTk2}\includegraphics[scale=0.5]{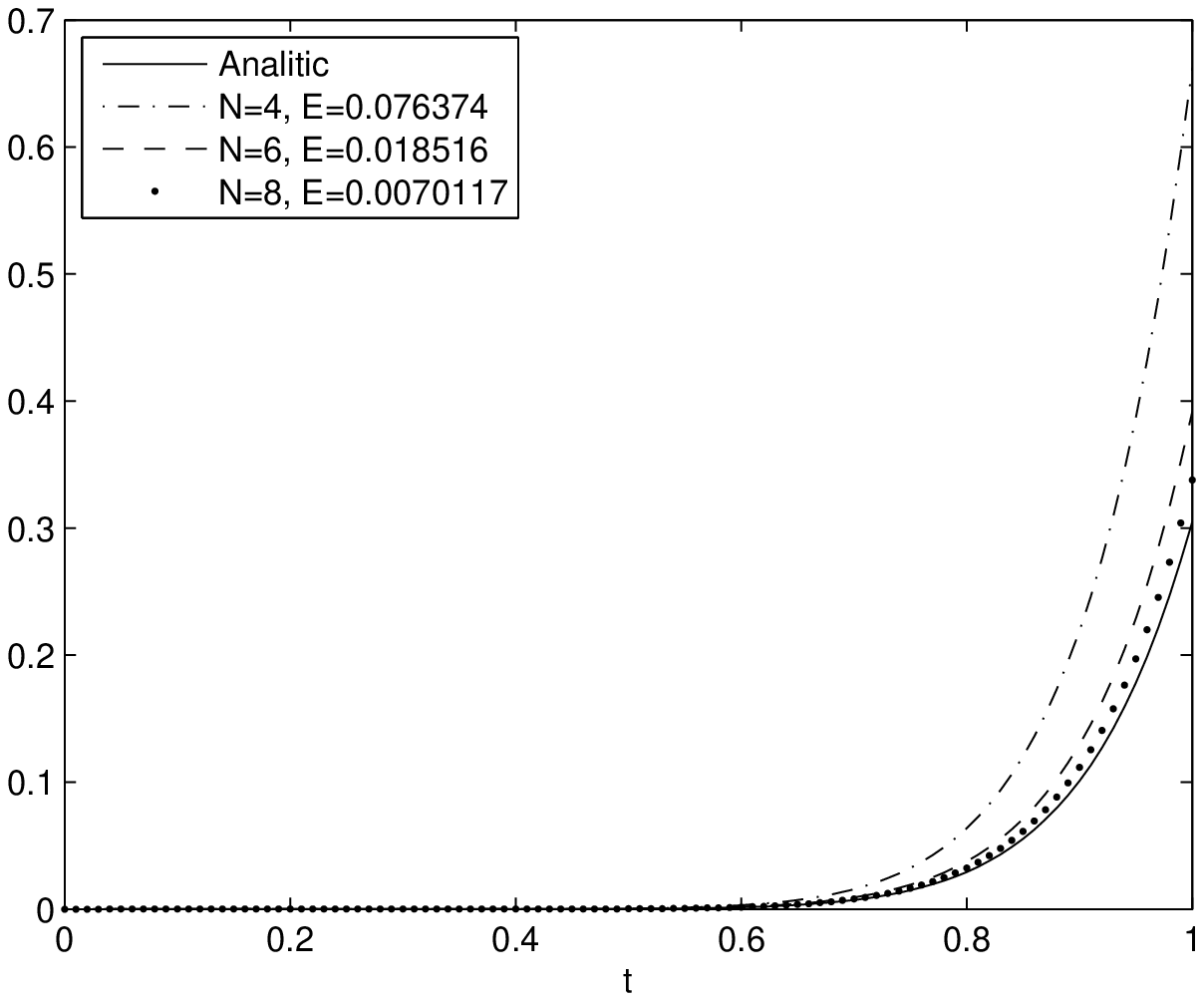}}
  \end{center}
  \caption{Analytic versus numerical approximation for a fixed $n$.}
  \label{Fig:IntExp1a}
\end{figure}

In this section we exemplify the proposed
approximation procedure with some examples.
In each step, we evaluate the accuracy of our method, \textrm{i.e.},
the error when substituting $\LI x$ by the approximation $\tilde{\LI}x$.
For that purpose, we take the distance given by
$$
E=\sqrt{\int_a^b \left(\LI x(t)-\tilde{\LI}x(t)\right)^2dt}.
$$
Firstly, consider $x_1(t)=t^3$ and $x_2(t)=t^{10}$ with $t\in[0,1]$. Then
$$
\LIz x_1(t)=\frac{\Gamma(4)}{\Gamma(4.5)}t^{3.5}\
\mbox{ and } \ \LIz x_2(t)=\frac{\Gamma(11)}{\Gamma(11.5)}t^{10.5}
$$
(\textrm{cf.} Property~2.1 in \cite{Kilbas}).
Let us consider Theorem~\ref{TheoDecomp} for $n=3$,
\textrm{i.e.}, expansion \eqref{def:n=3} for different values of step $N$.
For function $x_1$, small values of $N$ are enough
($N=3,4,5$). For $x_2$ we take $N=4,6,8$. In Figures~\ref{ExpTk}
and \ref{ExpTk2} we represent the graphs of the fractional integrals
of $x_1$ and $x_2$ of order $\a=0.5$ together with different
approximations. As expected, when $N$ increases
we obtain a better approximation for each fractional integral.


\begin{figure}[!ht]
  \begin{center}
    \subfigure[$\LIz(e^t)$]{\label{ExpEt}\includegraphics[scale=0.5]{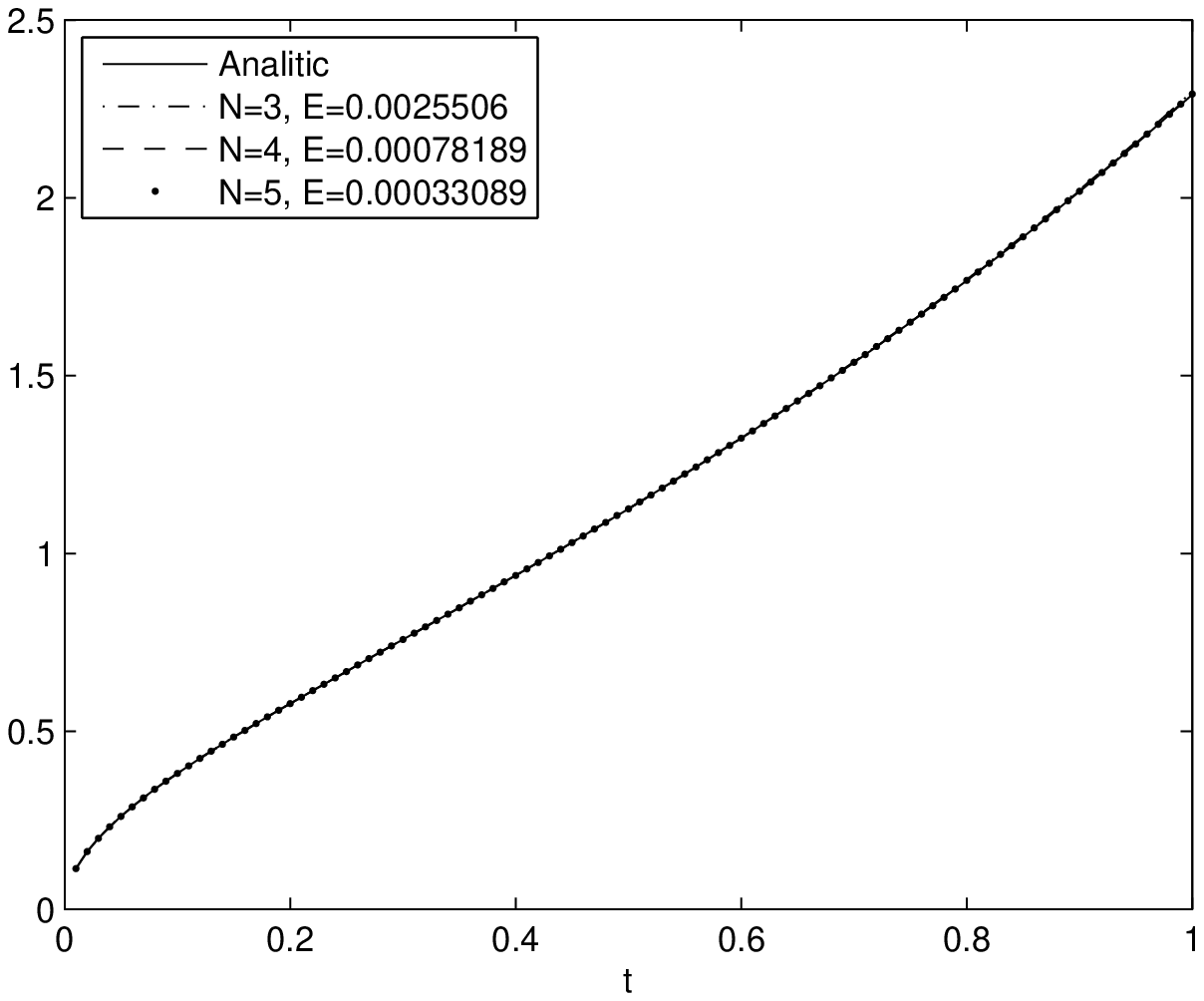}}
    \subfigure[$\LIz(\sin(t))$]{\label{ExpSint}\includegraphics[scale=0.5]{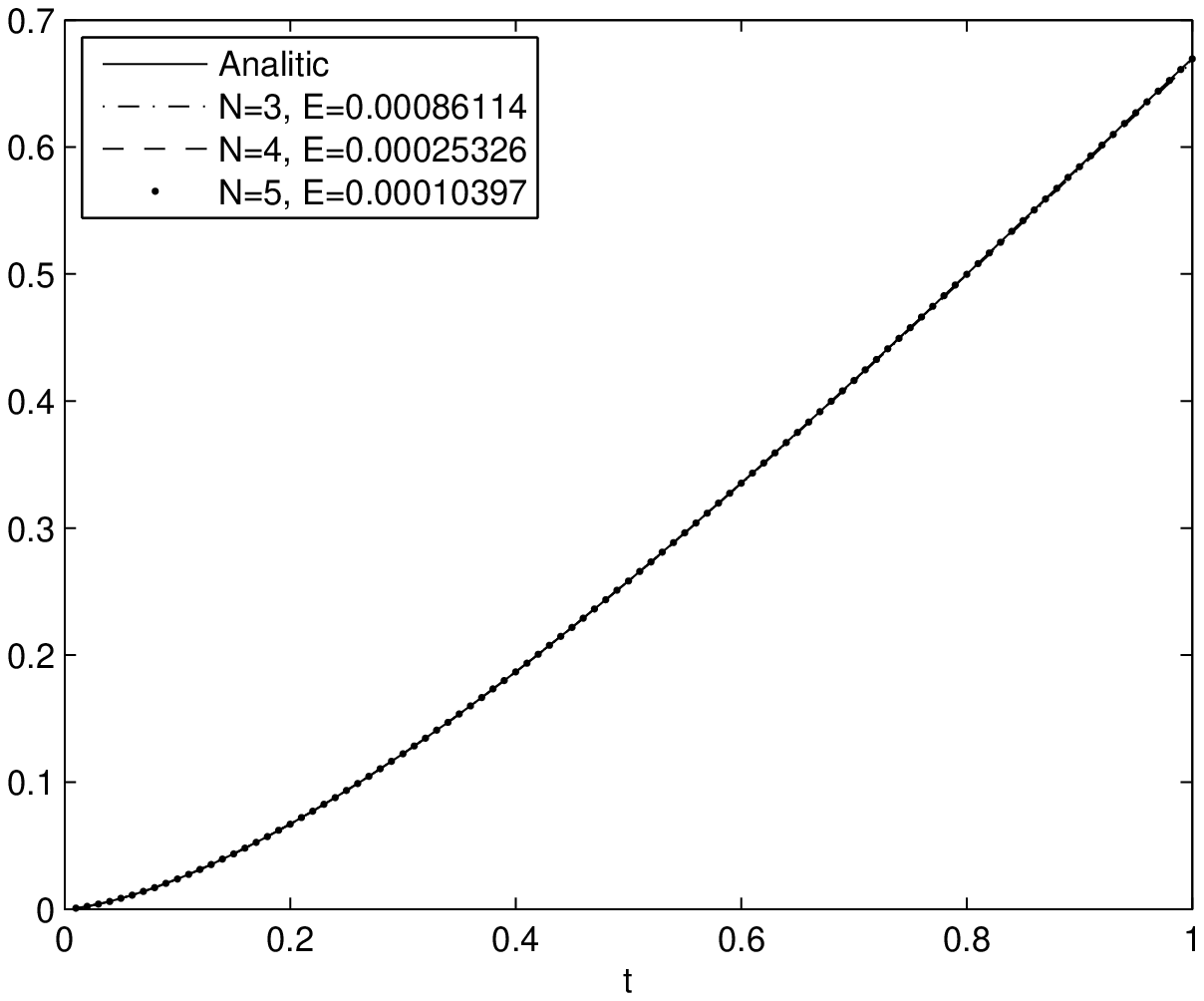}}
  \end{center}
  \caption{Analytic versus numerical approximation for a fixed $n$.}
  \label{Fig:IntExp2}
\end{figure}

Secondly, we apply our procedure to the transcendental functions
$x_3(t)=e^t$ and $x_4(t)=\sin(t)$. Simple calculations give
$$
\LIz x_3(t)=\sqrt{t}\sum_{k=0}^\infty\frac{t^k}{\Gamma(k+1.5)}
\ \mbox{ and }\  \LIz x_4(t)=\sqrt{t}\sum_{k=0}^\infty
\frac{(-1)^k t^{2k+1}}{\Gamma(2k+2.5)}.
$$
Figures~\ref{ExpEt} and \ref{ExpSint} show the numerical results
for each approximation, with $n=3$. We see that for a small value of $N$
one already obtains a good approximation for each function.


\begin{figure}[!ht]
\begin{center}
\subfigure[$\LIz(e^t)$]{\label{ExpEt2}\includegraphics[scale=0.5]{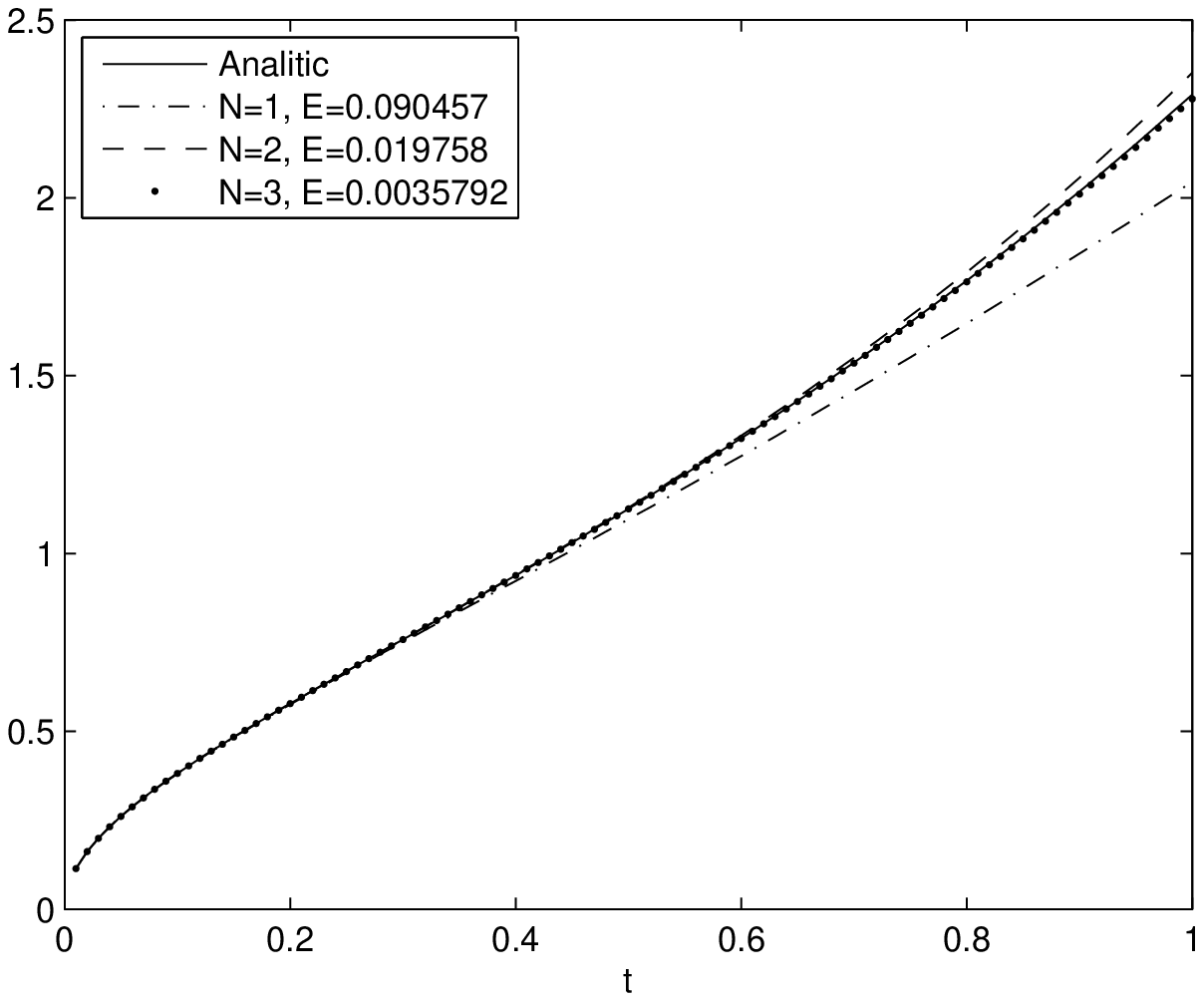}}
\subfigure[$\LIz(\sin(t))$]{\label{ExpSint2}\includegraphics[scale=0.5]{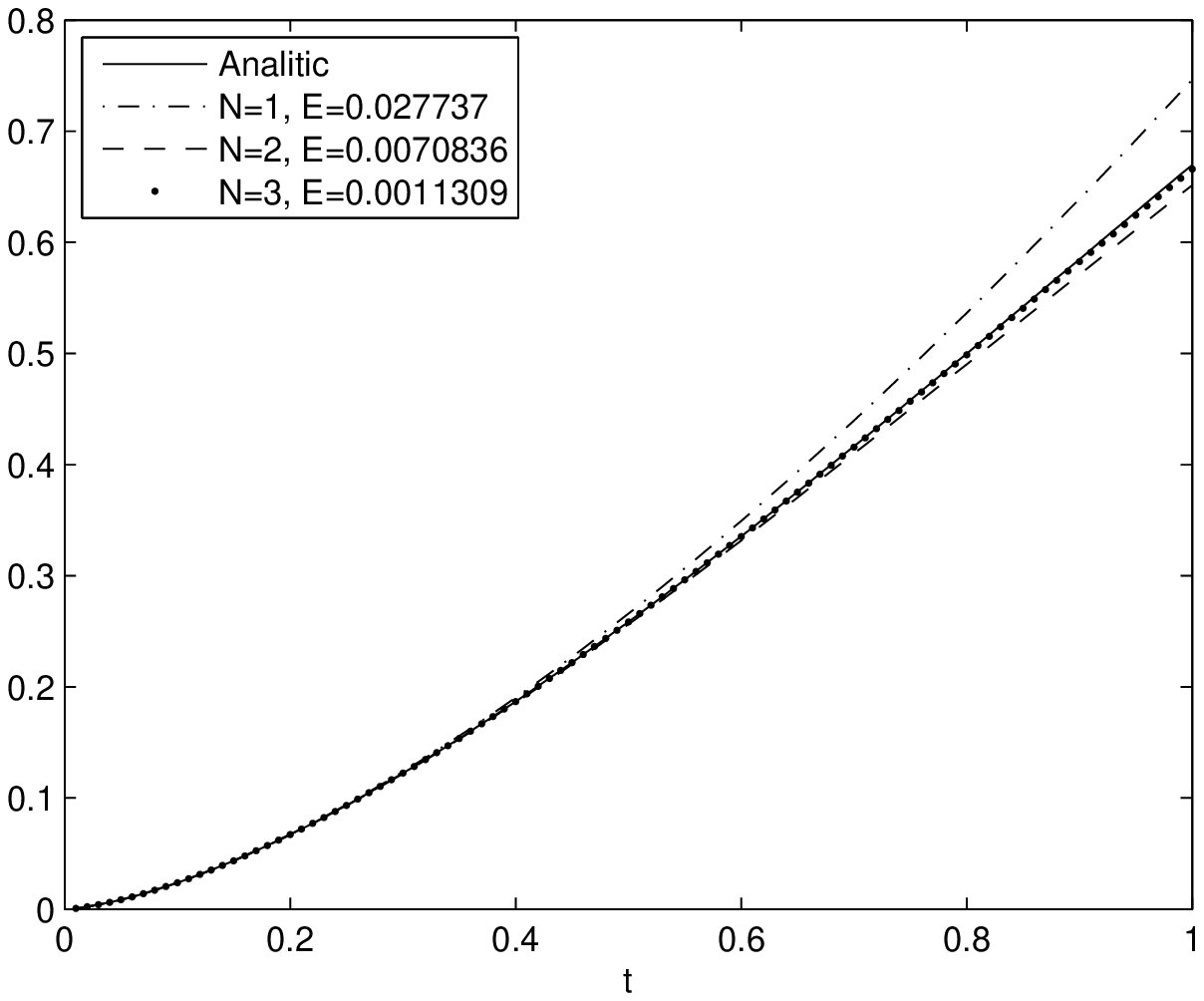}}
\end{center}
\caption{Numerical approximation using \eqref{analytical2} of previous literature.}
\label{Fig:IntExp3}
\end{figure}

For analytical functions, we may apply the well-known formula \eqref{analytical2}.
In Figure~\ref{Fig:IntExp3} we show the results of approximating
with \eqref{analytical2}, $N=1,2,3$, for functions $x_3(t)$
and $x_4(t)$. We remark that, when we consider expansions up to the second derivative,
\textrm{i.e.}, the cases $n=3$ as in \eqref{def:n=3} and expansion \eqref{analytical2}
with $N=2$, we obtain a better accuracy using our approximation
\eqref{def:n=3} even for a small value of $N$.


\begin{figure}[!ht]
\begin{center}
\subfigure[$\LIz(t^3)$]{\label{ExpTkn}\includegraphics[scale=0.5]{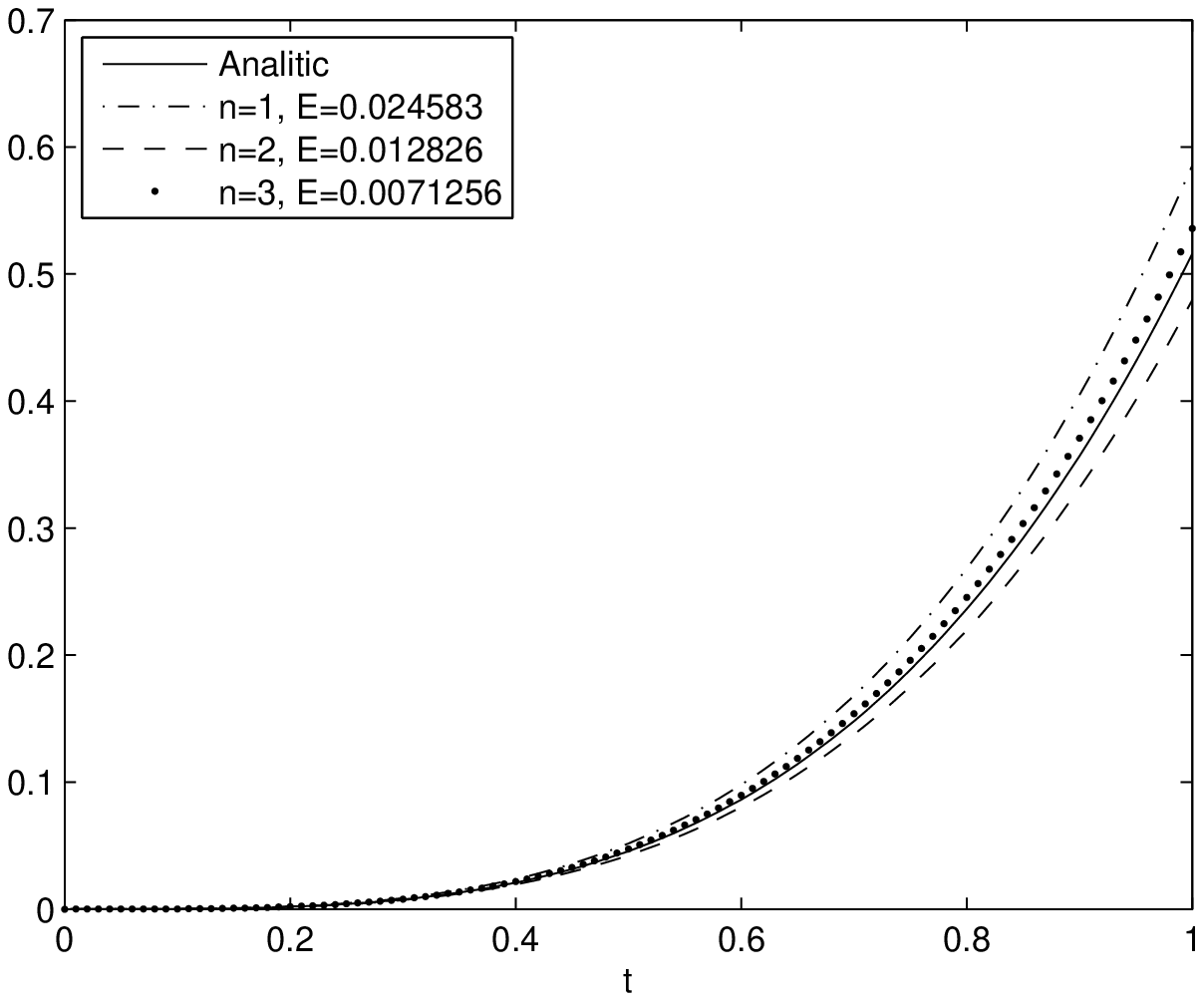}}
\subfigure[$\LIz(t^{10})$]{\label{ExpTkn2}\includegraphics[scale=0.5]{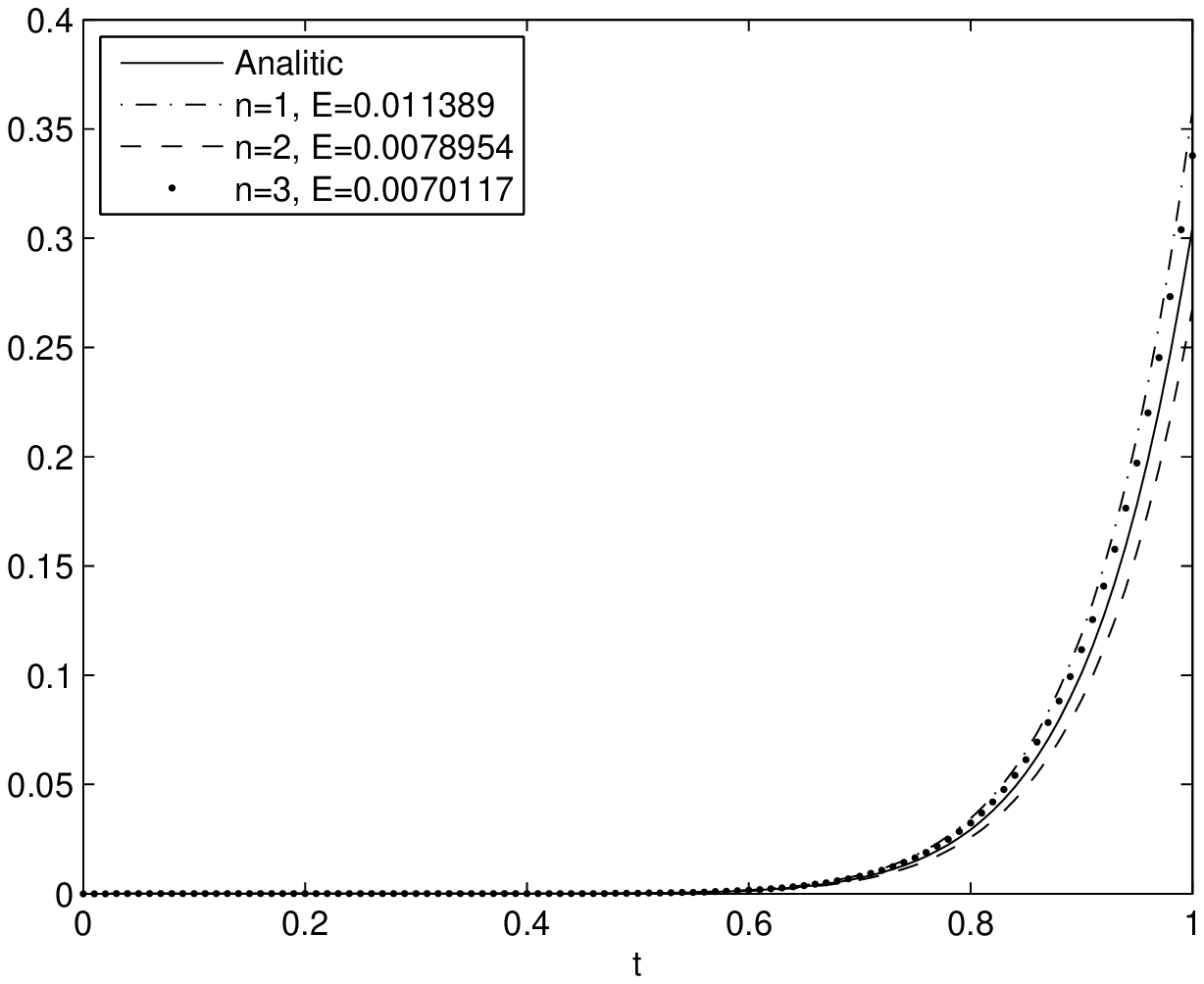}}
\end{center}
\caption{Analytic versus numerical approximation for a fixed $N$.}
\label{Fig:IntExp4}
\end{figure}

Another way to approximate fractional integrals is to fix $N$
and consider several sizes for the decomposition, \textrm{i.e.},
letting $n$ to vary. Let us consider the two test functions $x_1(t)=t^3$
and $x_2(t)=t^{10}$, with $t\in[0,1]$ as before. In both cases
we consider the first three approximations of the fractional integral,
\textrm{i.e.}, for $n=1,2,3$. For the first function we fix
$N=3$, for the second one we choose $N=8$. Figures~\ref{ExpTkn} and \ref{ExpTkn2}
show the numerical results. As expected, for a greater value of $n$ the error decreases.


\begin{figure}[!ht]
  \begin{center}
    \subfigure[$\LIz(t^3)$]{\label{ExpTkA}\includegraphics[scale=0.5]{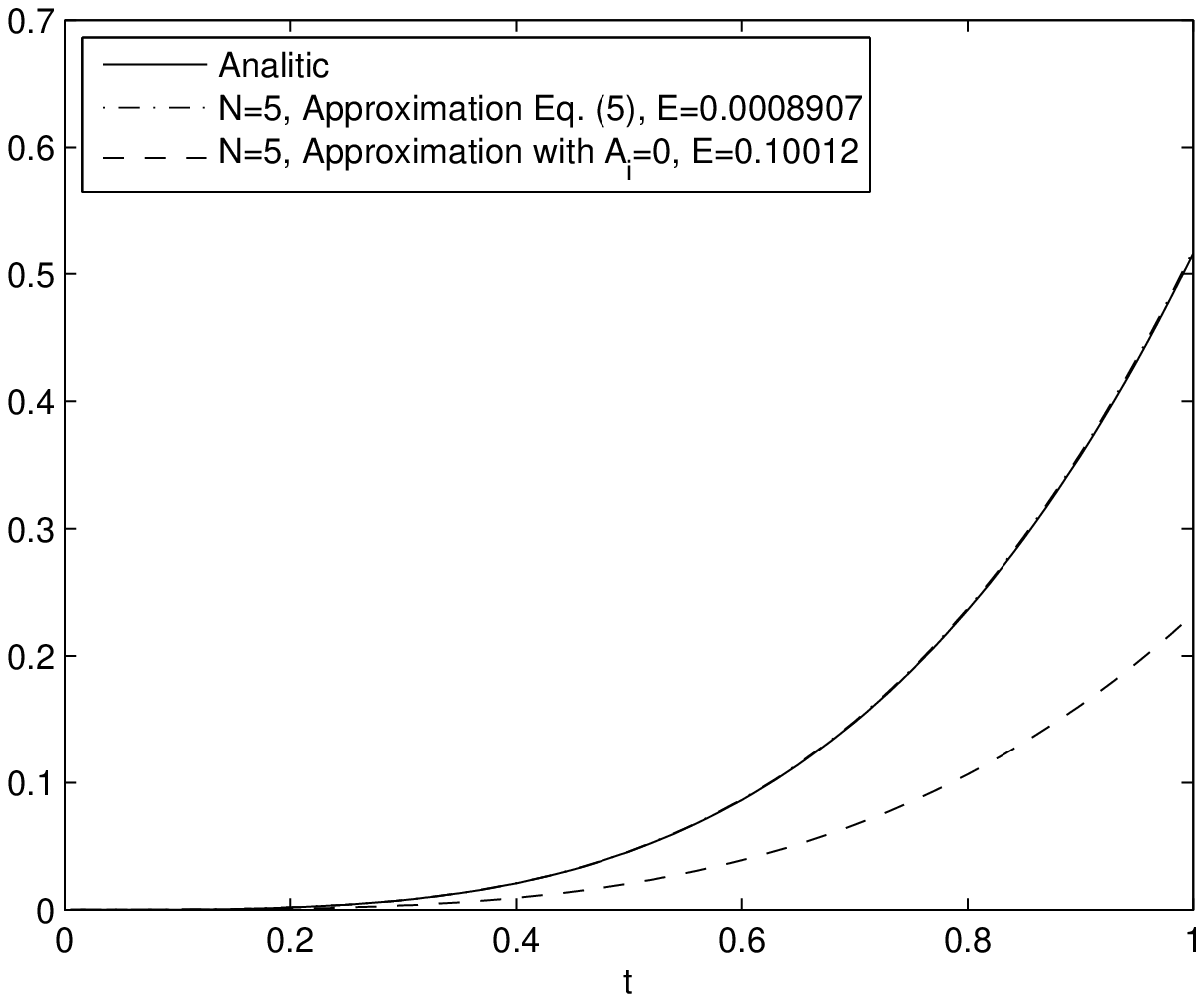}}
    \subfigure[$\LIz(t^{10})$]{\label{ExpTk2A}\includegraphics[scale=0.5]{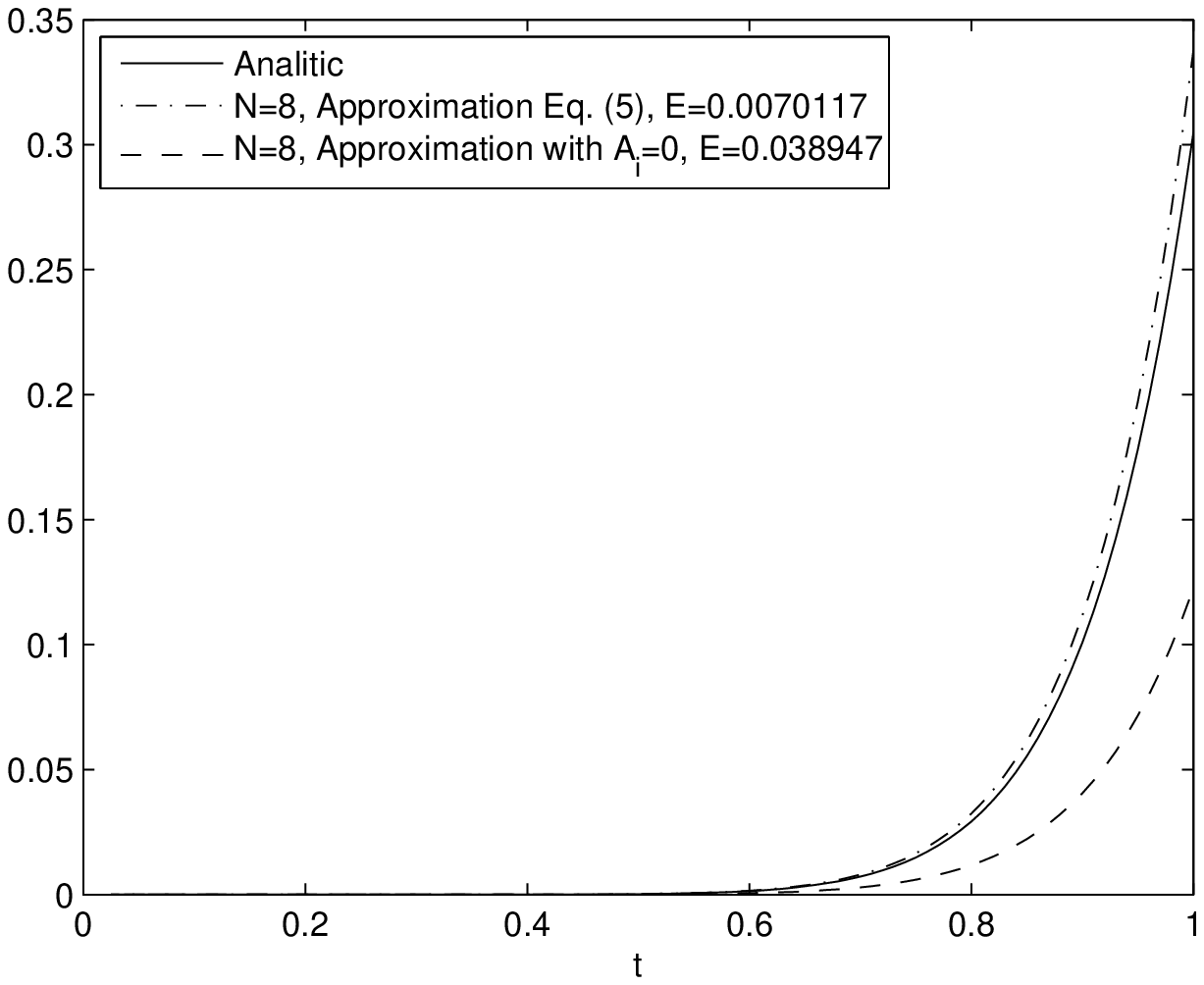}}
  \end{center}
  \caption{Comparison of approximation \eqref{def:n=3} and approximation with $A_i=0$.}
  \label{Fig:IntExp1b}
\end{figure}

We mentioned before that although the terms $A_i$ are all equal to zero,
for $i \in \{0,\ldots,n-1\}$, we consider them in the decomposition formula.
Indeed, after we truncate the sum, the error is lower. This is illustrated
in Figures~\ref{ExpTkA} and \ref{ExpTk2A}, where we study the approximations
for $\LIz x_1(t)$ and $\LIz x_2(t)$ with $A_i\not=0$ and $A_i=0$.


\section{Applications}
\label{secExp2}

In this section we show how the proposed approximations can be applied
into different subjects. For that, we consider a fractional integral
equation (Example~\ref{example:1}) and a fractional variational problem in which the
Lagrangian depends on the left Riemann--Liouville fractional integral
(Example~\ref{example:2}). The main idea is to rewrite the initial problem by
replacing the fractional integrals by an expansion
of type \eqref{analytical} or \eqref{ExpanDecomp},
and thus getting a problem involving integer derivatives,
which can be solved by standard techniques.


\begin{example}[Fractional integral equation]
\label{example:1}
Consider the following fractional system:
\begin{equation}
\label{system}
\left\{ \begin{array}{l}
\LIz x(t)=\frac{\Gamma(4.5)}{24}t^4\\
x(0)=0.
\end{array}\right.
\end{equation}
Since $\LIz t^{3.5}=\frac{\Gamma(4.5)}{24}t^\a$, the function
$t\mapsto t^{3.5}$ is a solution to problem \eqref{system}.

To provide a numerical method to solve such type of systems,
we replace the fractional integral by approximations \eqref{analytical2}
and \eqref{Approx:LI}, for a suitable order. We remark that the order
of approximation, $N$ in \eqref{analytical2} and $n$ in \eqref{Approx:LI},
are restricted by the number of given initial or boundary conditions.
Since \eqref{system} has one initial condition, in order to solve it numerically,
we will consider the expansion for the fractional integral up to the first derivative,
\textrm{i.e.}, $N=1$ in \eqref{analytical2} and $n=2$ in \eqref{Approx:LI}.
The order $N$ in \eqref{Approx:LI} can be freely chosen.

Applying approximation \eqref{analytical2}, with $\a=0.5$,
we transform \eqref{system} into the initial value problem
$$
\left\{
\begin{array}{l}
1.1285t^{0.5} x(t)-0.3761t^{1.5} x'(t)=\frac{\Gamma(4.5)}{24}t^4,\\
x(0)=0,
\end{array}\right.
$$
which is a first order ODE.
The solution is shown in Figure~\ref{IntFI}.
It reveals that the approximation remains close to the exact solution
for a short time and diverges drastically afterwards. Since we have
no extra information, we cannot increase the order of approximation to proceed.

To use expansion \eqref{ExpanDecomp}, we rewrite the problem as
a standard one, depending only on a derivative of first order.
The approximated system that we must solve is
$$
\left\{
\begin{array}{l}
A_0(0.5,N)t^{0.5} x(t)+A_1(0.5,N)t^{1.5} x'(t)
+\sum_{p=2}^NB(0.5,p)t^{1.5-p}V_p(t)
=\frac{\Gamma(4.5)}{24}t^4,\\
V_p'(t)=(p-1)t^{p-2}x(t),\quad p=2,3,\ldots,N,\\
x(0)=0,\\
V_p(0)=0, \quad p=2,3,\ldots,N,\\
\end{array}\right.
$$
where $A_0,A_1$ are given as in \eqref{Def:A} and $B$ is given
by Theorem~\ref{TheoDecomp}. Here, by increasing $N$, we get better
approximations to the fractional integral and we expect more accurate
solutions to the original problem \eqref{system}. For $N=2$ and $N=3$
we transform the resulting system of ordinary differential equations
to a second and a third order differential equation respectively. Finally,
we solve them using the \textsf{Maple} built in function \textsf{dsolve}.
For example, for $N=2$ the second order equation takes the form
$$
\left\{
\begin{array}{l}
V_2''(t)=\frac{6}{t}V_2'(t)+\frac{6}{t^2}V_2(t)-5.1542t^{2.5}\\
V(0)=0\\
V'(0)=x(0)=0
\end{array}
\right.
$$
and the solution is $x(t)=V_2'(t)=1.34t^{3.5}$.
In Figure~\ref{momFI} we compare the exact solution with numerical
approximations for two values of $N$.

\begin{figure}[!ht]
\begin{center}
\subfigure[Approximation by \eqref{analytical2}]{\label{IntFI}\includegraphics[scale=0.5]{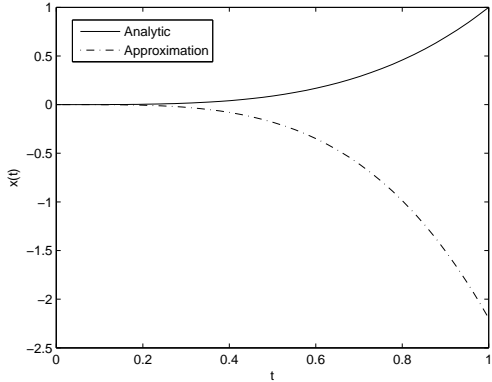}}
\subfigure[Approximation by \eqref{Approx:LI}]{\label{momFI}\includegraphics[scale=0.5]{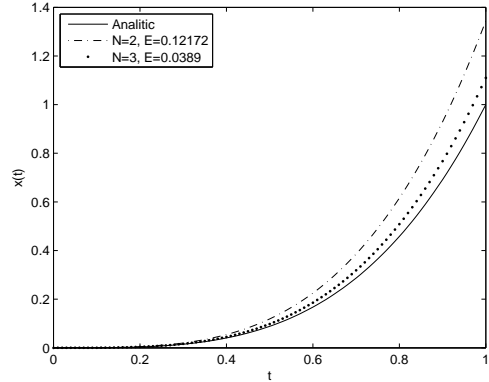}}
\end{center}
\caption{Analytic versus numerical solution to problem \eqref{system}.}
\label{Fig:IntExp5}
\end{figure}
\end{example}


\begin{example}[Fractional variational problem]
\label{example:2}
Let $\a\in (0,1)$. Consider the problem
\begin{equation}
\label{covEx}
\begin{gathered}
J[x(\cdot)]=\int_0^1 \left(\LIza x(t)-t\right)^2 dt \longrightarrow \min ,\\
x(0)=0, \quad x(1)=\frac{\Gamma(\a+1.5)}{\Gamma(1.5)}.
\end{gathered}
\end{equation}
Problems of the calculus of variations of type \eqref{covEx},
with a Lagrangian depending on fractional integrals,
were first introduced in \cite{MyID:148}.
In this case the solution is rather obvious. Indeed, for
\begin{equation}
\label{eq:sol:vp}
x(t)=\frac{\Gamma(1.5+\a)}{\Gamma(1.5)}\sqrt{t}
\end{equation}
one has $\LIza x(t)=t$ and $J[x(\cdot)] = 0$.
Since functional $J$ is non-negative,
\eqref{eq:sol:vp} is the global minimizer
of \eqref{covEx}.

Using \eqref{Approx:LI},
we approximate problem \eqref{covEx} by
\begin{equation}
\label{covApp}
\begin{gathered}
\tilde{J}[x(\cdot)]=\int_0^1 \left[A_0(\a,N)t^{\a} x(t)
+A_1(\a,N)t^{1+\a} x'(t)+\sum_{p=2}^N
B(\a,p)t^{1+\a-p}V_p(t)-t\right]dt \longrightarrow \min ,\\
V_p'(t)=(p-1)t^{p-2}x(t),\quad p=2,3,\ldots,N,\\
V_p(0)=0, \quad p=2,3,\ldots,N,\\
x(0)=0,\quad x(1)=\frac{\Gamma(\a+1.5)}{\Gamma(1.5)}.
\end{gathered}
\end{equation}
This is a classical variational problem,
constrained by a set of boundary conditions
and ordinary differential equations. One way to solve such
a problem is to reformulate it
as an optimal control problem \cite{Pontryagin}.
Let us introduce the control variable
$$
u(t) = A_0(\a,N)t^{\a} x(t)
+A_1(\a,N)t^{1+\a} x'(t)
+\sum_{p=2}^NB(\a,p)t^{1+\a-p}V_p(t).
$$
Then \eqref{covApp} becomes the classical optimal control problem
\begin{equation*}
\begin{gathered}
\tilde{J}[x(\cdot)]=\int_0^1 \left[u(t)-t\right]dt \longrightarrow \min ,\\
x'(t)=-A_0A_1^{-1}t^{-1}x(t)-\sum_{p=2}^N
A_1^{-1}B_pt^{-p}V_p(t)+A_1^{-1}t^{-1-\a}u(t),\\
V_p'(t)=(p-1)t^{p-2}x(t),\quad p=2,3,\ldots,N,\\
V_p(0)=0, \quad p=2,3,\ldots,N,\\
x(0)=0,\quad x(1)=\frac{\Gamma(\a+1.5)}{\Gamma(1.5)},
\end{gathered}
\end{equation*}
where $A_i=A_i(\a,N)$ and $B_p=B(\a,p)$.
For $\a=0.5$ and $N=2$, application of the
Hamiltonian system with multipliers
$\lambda_1(t)$ and $\lambda_2(t)$ \cite{Pontryagin},
gives the two point boundary value problem
\begin{equation}
\label{tpbvp}
\left\{
\begin{array}{rl}
x'(t)&=-A_0A_1^{-1}t^{-1}x(t)-A_1^{-1}B_2t^{-2}V_2(t)t^{1-\a}
-\frac{1}{2}A_1^{-2}t^{-2-2\a}\lambda_1(t)+A_1^{-1}t^{-\a},\\
V_2'(t)&=x(t),\\
\lambda_1'(t)&=A_0A_1^{-1}t^{-1}\lambda_1(t)-\lambda_2(t),\\
\lambda_2'(t)&=A_1^{-1}B_2t^{-2}\lambda_1(t),\\
\end{array}\right.
\end{equation}
with boundary conditions
\begin{equation*}
\label{sysB}
\left\{
      \begin{array}{l}
        x(0)=0, \\
        V_2(0)=0,
      \end{array}
\right.\qquad\left\{
      \begin{array}{l}
        x(1)=1,\\
        \lambda_2(1)=0.
      \end{array}
\right.
\end{equation*}
In practice the value of $n$ is chosen taking into account
the number of boundary conditions. For problem \eqref{covEx},
to avoid lack of boundary conditions, we take $n=2$.
On the other hand, the value of $N$ is only restricted by the computational power
and the efficiency of the numerical method chosen to solve the approximated problem.
Figure~\ref{Fig:IntExp6} shows the solution to system \eqref{tpbvp}
together with the exact solution to problem \eqref{covEx}
and corresponding values of $J$.


\begin{figure}[!ht]
\begin{center}
\includegraphics[scale=0.5]{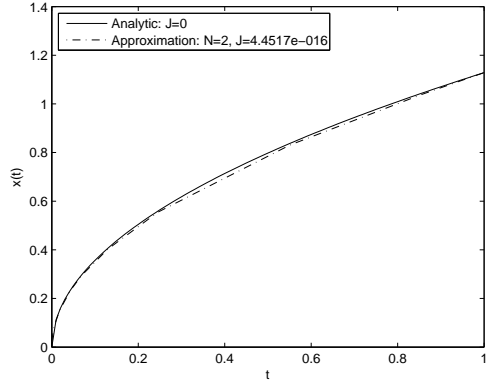}
\end{center}
\caption{Analytic versus numerical solution to fractional variational problem \eqref{covEx}.}
\label{Fig:IntExp6}
\end{figure}
\end{example}


\section*{Acknowledgments}

Work supported by {\it FEDER} funds through
{\it COMPETE} --- Operational Programme Factors of Competitiveness
(``Programa Operacional Factores de Competitividade'')
and by Portuguese funds through the
{\it Center for Research and Development
in Mathematics and Applications} (University of Aveiro)
and the Portuguese Foundation for Science and Technology
(``FCT --- Funda\c{c}\~{a}o para a Ci\^{e}ncia e a Tecnologia''),
within project PEst-C/MAT/UI4106/2011
with COMPETE number FCOMP-01-0124-FEDER-022690.
Pooseh was also supported by FCT through the Ph.D. fellowship
SFRH/BD/33761/2009.



\end{document}